\documentclass[12pt]{amsart}

\begin{document}
\title{Binomial theorem and exponent for variables commuting as $yx=qxy$}
\author{A. V. Stoyanovsky}
\address{Russian State University of Humanities}
\email{stoyan@mccme.ru}
\begin{abstract}
We state analogs of the binomial theorem and the exponential function for
variables $x$, $y$ commuting as $yx=qxy$.
\end{abstract}
\maketitle

Let $x$, $y$ be formal variables satisfying the commutation relation
$yx=qxy$, where $q$ is a number. In this note we are going to write out the formulas
for $(x+y)^n$ and for the function (formal series) $\exp_q(x)$ satisfying
the functional equation of the exponent. Note that the usual approach to $q$-analysis
is somewhat another (see, for example, [1,2]).

The author is grateful to V.~V.~Dolotin and A.~N.~Zhukov for stimulating questions.
Thanks are also due to Professors Tom Koornwinder and Uwe Franz, who informed me that this
result is not new, but has been discovered by Schutzenberger in 1953 [3], see also
[4,5] for further information and references.

\medskip

{\bf Theorem 1.}{\it
\begin{equation}
(x+y)^n=\sum_{k=0}^n{n\choose k}_qx^ky^{n-k},
\end{equation}
where ${n\choose k}_q$ are the usual $q$-binomial coefficients (Gaussian polynomials),
\begin{equation}
{n\choose k}_q=\frac{(n!)_q}{(k!)_q((n-k)!)_q},
\end{equation}
where
\begin{equation}
(n!)_q=(1)_q(2)_q(3)_q\ldots(n)_q,\ \  (k)_q=1+q+q^2+\ldots+q^{k-1},
\end{equation}
satisfying the recurrence relation
\begin{equation}
{n\choose k}_q={n-1\choose k}_q+q^{n-k}{n-1\choose k-1}_q,\ \ {n\choose 0}_q={n\choose n}_q=1,
\end{equation}
and the symmetry relation
\begin{equation}
{n\choose k}_q={n\choose n-k}_q.
\end{equation}
}
\bigskip

{\bf Theorem 2.} {\it If we put
\begin{equation}
\exp_q(x)=\sum_{n=0}^\infty\frac{x^n}{(n!)_q},
\end{equation}
then
\begin{equation}
\exp_q(x+y)=\exp_q(x)\exp_q(y).
\end{equation}
}

Proof of Theorem~1 is by induction. Theorem~2 follows from Theorem~1.

\end{document}